\documentclass[11pt,a4paper]{amsart}
\usepackage{amsmath,amssymb,latexsym,graphicx}
\newcommand{\bd}{{\mathbb{D}}} 
\newcommand{\bn}{{\mathbb{N}}}
\newcommand{\br}{{\mathbb{R}}}

\newcommand{\bc}{{\mathbb{C}}}

\newcommand{\bt}{{\mathbb{T}}}

\newcommand{\css}{{\mathcal{S}}}
\newcommand{\cv}{{\mathcal{V}}}

\newcommand{\s}{\sigma}

\renewcommand{\t}{\theta}
\renewcommand{\d}{\delta}

\newcommand{\g}{\gamma}

\newcommand{\z}{\zeta}

\newcommand{\nt}{\noindent}

\newcommand{\bsl}{\backslash}

\newcommand{\ti}{\tilde}

\newcommand{\lt}{\left}
\newcommand{\rt}{\right}

\DeclareMathOperator{\tr}{\rm tr}

\allowdisplaybreaks
\numberwithin{equation}{section}

\newtheorem{theorem}{Theorem}[section]
\newtheorem*{lemma}{Main Lemma}
\newtheorem{proposition}[theorem]{Proposition}

\newtheorem{conjecture}[theorem]{Conjecture}
\theoremstyle{definition}

\newtheorem{remark}[theorem]{Remark}

\begin{document}

\title[Absolutely continuous spectrum of a Schr\"odinger operator]
{Absolutely continuous spectrum of a Schr\"odinger operator on a tree} 
\author[S. Kupin]{S. Kupin}

\address{Universit\'e Aix-Marseille, 39, rue Joliot-Curie, 13453 
Marseille \\Cedex 13, France}
\email{kupin@cmi.univ-mrs.fr}

\date{May, 15,  2008}


\begin{abstract} 
We give sufficient conditions for the presence of the absolutely continuous 
spectrum of a Schr\"odinger operator on a regular rooted tree without loops 
(also called regular Bethe lattice or Cayley tree).
\end{abstract}

\maketitle

\vspace{-0.5cm}
\section*{Introduction and results} \label{s0}
The spectral properties of Schr\"odinger operators on graphs have numerous applications in physics and they have been intensively studied since late 90's.

We will be mainly interested in the properties of the absolutely continuous component of the spectral measure of  a discrete Schr\"odinger operator $H_V$ on a tree, see \eqref{e1} for an example. Probably, the first specific results in this direction were obtained by Klein \cite{kl1} who proved the presence of the absolutely continuous component   for $H_V$'s with random iid potential on a regular Bethe lattice. Recently, Aizenman-Sims-Warzel \cite{aiz1} obtained the result  with the help of a new general method. They also handled quasi-periodic operators \cite{aiz4} and a Laplacian on a random quantum tree \cite{aiz2}; see  Aizenman-Sims-Warzel \cite{aiz3} for a nice overview of the topic.  We also mention interesting papers by Froese-Hasler-Spitzer \cite{fro1, fro2} and Breuer \cite{br1}.

Almost simultaneously to the above-mentioned works,  Killip-Simon \cite{ks}, Nazarov-Peherstrorfer-Volberg-Yuditskii \cite{nvyu} obtained important results in the spectral theory of one-dimensional (1D) Schr\"odinger operators and, more generally, Jacobi matrices. These and subsequent papers \cite{ku1, ku2, ku3, sz1, zl1} gave a fairly complete picture of the spectral behavior of these 1D objects. 

It was hence very tempting to apply the well-developed methods of the one-dimensional analysis to the spectral problems for Schr\"odinger operators on trees. The first step in this direction was made by Denisov \cite{sd1}, who succeeded to carry over methods of Simon \cite{si1} to $H_V$ described in \eqref{e1}.

However,  the general picture remained quite unclear. In particular, we did not understand to what extent the construction for 1D worked for Schr\"odinger operators on trees. This gap is fixed by the present paper. Amongst other results, we prove the Main Lemma (see Section \ref{s1}) which expresses the Jost solutions for $H_V$ in terms of corresponding perturbation determinants. This observation implies immediately that there are strong parallels  between spectral behavior of 
1D Schr\"odinger operators and similar objects on trees, and we recover a big part of the 1D theory for these $H_V$'s. In particular, the sum rules of higher order for 1D Jacobi matrices become the ``sum inequalities"  for $H_V$. These higher order sum inequalities are proved in Theorems \ref{t2}, \ref{t3}.

\begin{figure}[t]
\begin{center}
\includegraphics[height=4cm]{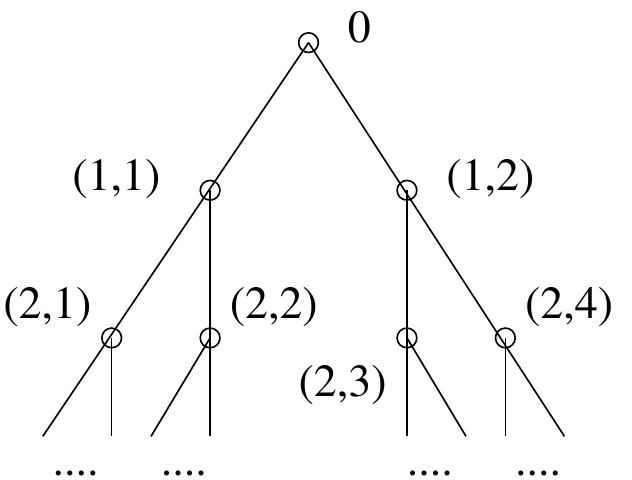}
\end{center}
\caption{The tree $T$}\label{f1}
\end{figure}

Let $T=T_0$ be a regular binary tree of without loops (also called Cayley tree or Bethe lattice). Its root vertex is denoted by 0.  The set of all vertices of the tree is denoted by $\cv(T)$. The distance $|x_1-x_2|$ 
between 
two vertices $x_1,x_2\in\cv(T)$ is the number of edges of the (unique) path 
leading from $x_1$ to $x_2$.  The sphere of radius $n$ and centered at 0 is
$$
\css(0,n)=\{x\in\cv(T): |x|=|x-0|=n\}.
$$
Every vertex in the tree has one ascendant and two descendants. The descendants 
of different vertices are different since the tree does not have closed loops. 
So, the sphere $\css(0,n)$ contains $2^n$ vertices $x=(n,k),\ k=1,\dots, 2^n$. 
That is,
\begin{itemize}
\item $\cv(T)=\{0,(n,k): n\in\bn,\ k=1,\dots, 2^n\}$,\
\item the descendants of $0$ are vertices $(1,1), (1,2)$,
\item for $n\ge 1$, the descendants of $(n,j), j=1,\dots, 2^n$, are $(n+1, 
2j-1)$ and $(n+1,2j)$, see Figure \ref{f1}.
\end{itemize}  

Let 
\begin{equation}\label{e01}
l^p(T)=\{\{u(x)\}_{x\in\cv(T)}:\ \sum_{x\in\cv(T)} |u(x)|^p<\infty\}
\end{equation}
with $1\le p\le\infty$. The standard ``basis'' vectors are $\{e_x\}_{x\in\cv(T)}$, 
where 
$e_x(x)=1$ and $e_x(y)=0$ for $y\not=x,\ y\in\cv(T)$. 
The free Laplacian $H_0=H_{0,T}$ is defined as
$$
(H_0f)(x)=\sum_{x':|x'-x|=1} f(x'),
$$
where $f\in l^2(T)$. We also set 
$$
m_{H_0}(z)=((H_0-z)^{-1} e_0,e_0)=\int_\br \frac{d\mu_0}{x-z}
$$
to be the Weyl-Titchmarsh function of the operator. The Borel measure $\mu_0$ is 
called the spectral measure of $H_0$ (with respect to $e_0$). The spectrum 
$\s(H_0)$ coincides with $\mathrm{supp}\, \mu_0$, and 
$\s(H_0)=\s_{ac}(H_0)=[-2\sqrt 2,2\sqrt 2]$, see for example \cite[Sect.~2]{sd1}.

A Schr\"odinger operator on $l^2(T)$ is a diagonal perturbation of $H_0$,
\begin{equation}\label{e1}
(H_V f)(x)=\sum_{x':|x'-x|=1} f(x')+V(x).
\end{equation}
We always assume that $V=\{V(x)\}_{x\in\cv(T)}$ lies in $c_0(T)$, where
$$
c_0(T)=\{\{u(x)\}_{x\in\cv(T)}:\ \lim_{|x|\to+\infty} u(x)=0\}.
$$ 
Then the operator $H_V$ is self-adjoint and, once again, we define its spectral measure 
$\mu=\mu_{H_V}$ as
$$
m_{H_V}(z)=((H_V-z)^{-1} e_0,e_0)=\int_\br \frac{d\mu}{x-z}.
$$
Since $H_V-H_0$ is compact, the Weyl-von Neumann theorem says that the essential spectrum
$\s_{ess}(H_V)$ of the operator $H_V$ equals $[-2\sqrt 2, 2\sqrt 2]$, and the point spectrum $\s_p(H_V)\subset \br\bsl \s_{ess}(H_V) $ accumulates to the points $\pm 2\sqrt 2$ only. It is convenient to enumerate $\s_p(H_V)=\{x^\pm_{V,T; s}\}$ as follows
\begin{equation}\label{e11}
x^-_{V, T; 1} \le \ldots \le x^-_{V, T; s}\le \dots <-2\sqrt 2,
\end{equation}
and 
\begin{equation}\label{e12}
2\sqrt 2<\dots \le x^+_{V, T; s}\le \dots \le x^+_{V, T; 1}.
\end{equation}
The numbering takes into account the (geometric) multiplicities of the eigenvalues.

For a given potential $V=\{V(x)\}_{x\in\cv(T)}$, define its ``truncation'' as
$$
V(n)=\{V(n;x)\}_{x\in\cv(T)}=
\left\{
\begin{array}{lcr}
V(x),&& |x|\le n,\\
0, && |x|>n.
\end{array}\right.
$$
Let $T_x, x\in\cv(T)$, be a subtree of $T$ growing from the vertex $x$. By $H_{V,T_x}$ 
we mean the Schr\"odinger operator with potential $V\big |_{T_x}$, the restriction of 
the original potential $V$ to $T_x$.  The notation $\s_p(H_{V,T_x})=\{x^\pm_{V,T_x; s}\}$ are self-obvious and stay for the point spectrum and eigenvalues of the operator $H_{V,T_x}$.

We give sufficient conditions for the support $\s_{ac}(H_V)$ of the absolutely 
continuous part of the measure $\mu_{H_V}$ to fill in the interval $[-2\sqrt 
2,2\sqrt 2]$. For instance, the following theorem is proved in Denisov \cite{sd1}. 

\begin{theorem}\label{t1} Let $H_V$ be a Schr\"odinger operator \eqref{e1}, $V\in c_0(T)$, and
$$
\sum^\infty_{n=1} \frac1{2^n}\sum_{x: |x|=n} V(x)^2<\infty.
$$
Then
\begin{eqnarray}
&& \int^{2\sqrt 2}_{-2\sqrt 2} \log\mu'(x)\cdot \sqrt{8-x^2}\, 
dx>-\infty,\nonumber\\
&& \limsup_n \Big (EV^{3/2}_{V(n),T}-\sum^\infty_{k=1} \frac 1{2^k}\sum_{x: |x|=k} 
EV^{3/2}_{V(n), T_x}\Big) <\infty, \label{e2}
\end{eqnarray}
where, for $p\ge 1$,
\begin{equation}\label{e3}
EV^p_{V(n),T}=\sum_s |x^+_{V(n),T; s}-2\sqrt 2|^p+ \sum_s |x^-_{V(n),T; s}+2\sqrt 
2|^p.
\end{equation}
\end{theorem}  
Above, $\mu'$ is the density of the absolutely continuous part of the measure 
$\mu$. Notice that the expression at the LHS of \eqref{e2} is actually non-negative. 

For a given $V\in l^\infty(T)$, define  $\d V$ as
$$
(\d V)(n,j)=
\lt\{
\begin{array}{lcl}
V(n-1,i)-V(n,2i),&& j=2i,\\
V(n-1,i)-V(n,2i-1),&& j=2i-1,
\end{array}
\rt.
$$
 
We prove the following theorems.
\begin{theorem}\label{t2} Let $H_V$ be a Schr\"odinger operator \eqref{e1}, $V\in c_0(T)$, and
\begin{eqnarray}
&& \sum^\infty_{n=1} \frac1{2^n}\sum_{x: |x|=n} V(x)^4<\infty,\quad 
 \sum^\infty_{n=2} \frac1{2^n}\sum_{x: |x|=n} (\d V)(x)^2<\infty, \nonumber\\
&&  \limsup_n \Big(EV^{5/2}_{V(n),T}-\sum^\infty_{k=1} \frac 1{2^k}\sum_{x: |x|=k} 
EV^{5/2}_{V(n), T_x}\Big) <\infty. \label{e31}
\end{eqnarray}
Then 
\begin{equation}\label{e32}
\int^{2\sqrt 2}_{-2\sqrt 2} \log\mu'(x)\cdot (8-x^2)^{3/2}\, dx>-\infty.
\end{equation}
\end{theorem}

\begin{theorem}\label{t3} Let $H_V$ be as in \eqref{e1}, $V\in c_0(T)$, and 
\begin{equation}\label{e35}
\sum^\infty_{n=1} \frac1{2^n}\sum_{x: |x|=n} V(x)^6<\infty,\quad  
\sum^\infty_{n=1} \frac1{2^n}\sum_{x: |x|=n} (\d V)(x)^2<\infty.
\end{equation}
Then we have
\begin{eqnarray}
&& \int^{2\sqrt 2}_{-2\sqrt 2} \log\mu'(x)\cdot (8-x^2)^{5/2}\, dx>-\infty, \label{e36}\\
&& \limsup_n \Big(EV^{7/2}_{V(n),T}-\sum^\infty_{k=1} \frac 1{2^k}\sum_{x: |x|=k} 
EV^{7/2}_{V(n), T_x}\Big)<\infty. \label{e4}
\end{eqnarray}
\end{theorem}
These theorems lead to conjectures stated in Section \ref{s2}.

\begin{remark}\label{r1} As in \cite{sd1}, the above theorem can be modified to assert that that relations \eqref{e35} (with $V\in l^\infty(T)$) yield \eqref{e36}, and, consequently, $[-2\sqrt 2, 2\sqrt 2]\subset \s_{ac}(H_V)$. The same applies to Conjecture \ref{c1}.p with odd $p$'s.
\end{remark}

Using results of Borichev-Golinskii-Kupin \cite{bgk}, we can express relations 
\eqref{e2}, \eqref{e31} and \eqref{e4} in terms of $\s_p(H_V)$ and the point spectra 
$\s_p(H_{V,T_x})$ of the corresponding operators.  Of course, the assumptions on the 
potential $V$ become considerably more stringent.
\begin{proposition}\label{p1} For $p\ge 1$, let $V\in l^{p'}(T), \ p'<p+1/2$.  Then the 
limit below exists and
$$
\lim_{n\to\infty} \Big\{EV^{p+1/2}_{V(n),T}-\sum^\infty_{k=1} \frac 1{2^k}\sum_{x: 
|x|=k} EV^{p+1/2}_{V(n), T_x}\Big\}=EV^{p+1/2}_{V,T}-\sum^\infty_{k=1} \frac 
1{2^k}\sum_{x: |x|=k} EV^{p+1/2}_{V, T_x}.
$$
We also have
$$
EV^{p+1/2}_{V,T}\le C(p,p',||V||_\infty) ||V||^{p'}_{p'},
$$
where $||.||_\infty, ||.||_p$ are the norms of $l^\infty(T), l^p(T)$, respectively.
\end{proposition} 

\section{Sketch of the proof of Theorem \ref{t2}}\label{s1}
As usual, we uniformize the domain $\bar\bc\bsl[-2\sqrt 2, 2\sqrt 2]$ with the help of 
the maps $z(\z)=\sqrt 2(\z+1/\z),\ \z(z)=\frac 1{2\sqrt 2}(z-\sqrt{z^2-8})$, where 
$\z\in\bd=\{\z: |\z|<1\}$.

Let for the moment $rank\, V<\infty$. We put for an arbitrary subtree $T_x\subset T, 
x\in\cv(T),$
$$
L_{T_x}(\z)=L_{T_x}(z(\z))=\det(H_{V, T_x}-z)(H_{0, T_x}-z)^{-1},
$$
where $H_{0, T_x}$ is the free Laplacian on $T_x$. Furthermore, let $\ti X$ be a 
finite subset of vertices of the tree $T$ with the property $T_x\cap T_y=\emptyset$ 
for $x,y\in\ti X$ and $x\not=y$. Obviously, we can speak about $H_{T_{\ti 
X}}=\oplus_{x\in\ti X} H_{T_x}$ and
$$
L_{T_{\ti X}}(\z)=L_{T_{\ti X}}(z(\z))=\det(H_{V, T_{\ti X}}-z)(H_{0, T_{\ti 
X}}-z)^{-1}=
\prod_{x\in\ti X} L_{T_x}(\z).
$$  
Consider the path $\g_y$ leading from 0 to $y\in\cv(T)$ and denote by $\ti X(y)$ the 
set of vertices lying on the distance one from vertices of the path, that is,
$$
\ti X(y)=\{x\in\cv(T): \exists w\in\cv(T)\cap\g_y, |x-w|=1\}.
$$
It is easy that $\ti X(y)$ has the above-mentioned disjointness property and, 
moreover,
$(\cv(T)\cap\g_y)\cup \cv(T_{\ti X(y)})=\cv(T)$.

The next lemma is the key to the proofs of Theorems \ref{t1}-\ref{t3}. 
It is new and it expresses the Jost solution of the operator $H_V$ in terms of 
$L_{T_{\ti X(y)}}$, compare to \cite[Theorem 2.16]{ks}. It goes without saying that the lemma holds 
also for ``sparse" trees considered in \cite{br1}.

\begin{lemma}\label{l1}
Let $rank\, V<\infty$ and $H_V$ be the Schr\"odinger operator \eqref{e1}. Let 
$f(\z)=\{f_y(\z)\}_{y\in T}\in l^2(T)$ and $f=(H_V-z(\z))^{-1}e_0$. Then, for $n=|y|$,
\begin{equation*}
f_y(\z)=\lt(\frac\z{\sqrt 2}\rt)^n\, L_{T_{\ti X(y)}}(\z)/L_T(\z).
\end{equation*}
\end{lemma}

The proofs of the theorem use the techniques developed in \cite{ks,nvyu,ku1,ku2} and 
the lemma.

\medskip\nt
{\it Sketch of the proof of Theorem \ref{t2}.}\ 
Let the potential $V\in c_0(T)$ satisfy the assumptions of the theorem.  
We do the computations for the operator $H_{V(N)}$, and then pass to the limit with 
respect to $N\to\infty$.

Make  the change of variables $z(\z)=\sqrt 2(\z+\z^{-1})$ and transfer the spectral 
measure $\mu_N=\mu_{H_{V(N)}}$ to the unit disk $\bd$ and its boundary. The absolutely continuous part of the image of the measure is then supported on the unit circle and its density is still 
denoted $\mu'_N$. We write $\{\z_{V(N),T_x;\, s}\}_s$ for the images of $\{x^\pm_{V(N),T_x;\, s}\}_s$. Then relations \eqref{e31}, \eqref{e32} read as 
\begin{equation*}
\int^{2\pi}_0 \log\mu'_N(e^{i\t})\, \sin^4\t\, d\t>-\infty,\ \ 
\limsup_N \Big(F_{V(N),T}-\sum^\infty_{k=1} \frac 1{2^k}\sum_{x: |x|=k} F_{V(N), 
T_x}\Big)<\infty, 
\end{equation*}
where
$$
F_{V(N),T}=\sum_s (1-|\z_{V(N),T_x; s}|)^5.
$$
 
For a given vertex $x\in\cv(T)$, we consider $H_{V(N), T_x}$ and the perturbation 
determinant $L_{T_x}(\z)$ ($=\det(H_{V(N), T_x}-z)(H_{0, T_x}-z)^{-1}$). 
The eigenvalues $\{\z_{V(N), T_x; s}\}$ coincide with the zeros of the determinant up to multiplicities. We have in a neighborhood of 
$\z=0$
$$
\log L_{T_x}(\z)=-\sum^\infty_{k=1} \frac1k \tr \lt(T_k(\frac 1{\sqrt 2} H_{V(N), 
T_x})-
T_k(\frac 1{\sqrt 2} H_{0, T_x})\rt)\, \z^k,
$$
where $T_k(2\cos\t)=2\cos k\t,\ k=0,1,2,\dots, $ are properly normalized Chebyshev 
polynomials of the first kind. The following identities hold: for $n=0$,
$$
\frac1{4\pi}\int^{2\pi}_0\log |L_{T_x}(e^{i\t})|^2\, d\t=\sum_s \log 1/|\z_{V(N), T_x; 
s}|,
$$
for $n\ge 1$,
\begin{eqnarray*}
\frac1{2\pi}\int^{2\pi}_0\log |L_{T_x}(e^{i\t})|^2\cos n\t\, d\t&=&\frac1n 
\sum_s \lt(1/|\z_{V(N), T_x; s}|^n-|\z_{V(N), T_x; s}|^n\rt)\\
&-&\frac1n \tr \lt(T_n(\frac 1{\sqrt 2} H_{V(N), T_x})-T_n(\frac 1{\sqrt 2} H_{0, 
T_x})\rt).
\end{eqnarray*}
Combining these equalities, we get
\begin{eqnarray}\label{e5}
&&\\
&&\frac1{2\pi}\int^{2\pi}_0\log |L_{T_x}(e^{i\t})|^2\, (16\sin^4\t)\, d\t=
\sum_s G(|\z_{V(N), T_x; s}|)\nonumber\\
&-&\frac18\tr\Big\{(H_{V(N), T_x}^4-24 H_{V(N), T_x}^2)-(H_{0, T_x}^4-24 H_{0, 
T_x}^2\Big\},\nonumber
\end{eqnarray}
where
$$
G(|\z|)=\frac 12\lt\{\lt(\frac 1{|\z|^4}-|\z|^4\rt)-8\lt(\frac 
1{|\z|^2}-|\z|^2\rt)+24\log |\z|\rt\}.
$$
This relation readily implies that
$$
\sum_s G(|\z_{V(N), T_x; s}|)\asymp\sum_s (1-|\z_{V(N), T_x; s}|)^5.
$$
Turning back to the operator $H_{V(N)}$ and its spectral characteristics, we observe 
that
$$
M_{H_{V(N)}}(\z)=-m_{H_{V(N)}}(z(\z))=f_0(\z),
$$
and the computation for $\mathrm{Im}\, M_{H_{V(N)}}(\z),\ \z=e^{i\t}\in\bt,$ gives
\begin{equation*}
\mathrm{Im}\, M_{H_{V(N)}}(\z)=\sqrt2\sin\t\frac{\sum_{y:|y|=N} |L_{T_{\ti X(y)}}(\z)|^2}{2^N  |L_T(\z)|^2}.
\end{equation*}
The inequality between the arithmetic and the geometric mean 
($\frac1n(a_1+\dots+a_n)\ge (a_1a_2\dots a_n)^{1/n}$ with $a_j\ge0$) and some simple 
combinatorics yield
\begin{eqnarray*}
\log\frac{\mathrm{Im}\, M_{H_{V(N)}}(\z)}{\sqrt 2\sin\t}&\ge&\frac1{2^N}\, 
\sum_{y:|y|=N} \log |L_{T_{\ti X(y)}}(\z)|^2-\log |L_T(\z)|^2\\
&=&\sum^N_{j=1}\frac 1{2^j} \sum_{x:|x|=j} \log|L_{T_x}(\z)|^2-\log |L_T(\z)|^2.
\end{eqnarray*}
We now apply equality \eqref{e5} to the logarithms in the RHS and transfer the sums corresponding to the point spectra to the LHS of the inequality. So we come to 
\begin{eqnarray}\label{e6}
&&\\
&&\frac1{2\pi}\int^{2\pi}_0\log \frac{\mathrm{Im}\, M_{H_{V(N)}}(e^{i\t})}{\sqrt 2\sin\t }\, 
(16\sin^4\t)\, d\t \nonumber\\
&+&\Big\{\sum_s G(|\z_{V(N),T_x; s}|)-\sum^N_{j=1}\frac 1{2^j} \sum_{x:|x|=j} \sum_s G(|\z_{V(N),T_x; s}|)\Big\}\nonumber\\
&\ge&\frac 18\tr\Big[K(H_{V(N), T})-K(H_{0,T})- \sum^N_{j=1}\frac 1{2^j} \sum_{x:|x|=j}
(K(H_{V(N), T_x})-K(H_{0,T_x})\Big], \nonumber
\end{eqnarray}
where $K(H)=H^4-24H^2$. An elementary calculation shows that the expression in the RHS of the relation satisfies the inequality
$$
[\ \dots ] \ge -C\lt\{\sum^N_{j=1}\frac 1{2^j} \sum_{x:|x|=j} V^4(x)+
\sum^N_{j=2}\frac 1{2^j} \sum_{x:|x|=j} (\d V)(x)^2\rt\}.
$$
Now,  take $\limsup_N$ of the both sides of inequality \eqref{e6}. Its RHS is finite by the assumptions of the theorem. Use the semi-continuity of the entropy (see \cite{ks}) to get
$$
\limsup_N \int^{2\pi}_0\log \frac{\mu'_N(e^{i\t})}{\sqrt 2\sin\t}\, \sin^4\t\, d\t\le
\int^{2\pi}_0\log \frac{\mu'(e^{i\t})}{\sqrt 2\sin\t }\, \sin^4\t\, d\t.  
$$ 
Hypothesis \eqref{e31} of the theorem says that
$$
\limsup_N \Big\{\sum_s G(|\z_{V(N),T_x; s}|)-\sum^N_{j=1}\frac 1{2^j} \sum_{x:|x|=j} \sum_s G(|\z_{V(N),T_x; s}|)\Big\}<\infty.
$$ 
The proof is complete. \hfill $\Box$

\section{Some open questions and conjectures}\label{s2}
The following conjecture seems very natural.
\begin{conjecture}\label{c1} \hspace{-1.5mm}{\bf p.}
Let $H_V$ be a Schr\"odinger operator  \eqref{e1} on a tree and $V\in c_0(T)$. Let, for an odd $p\ge 1$,
\begin{equation}\label{e7}
\sum^\infty_{n=1}\frac 1{2^n}\sum_{x:|x|=n} V(x)^{2p}<\infty, \quad
\sum^\infty_{n=2}\frac 1{2^n}\sum_{x:|x|=n} (\d V)(x)^2<\infty,
\end{equation}
Then
\begin{eqnarray}
&& \int^{2\sqrt 2}_{-2\sqrt 2} \log\mu'(x)\, (8-x^2)^{p-1/2}\, dx>-\infty, \nonumber\\
&&\limsup_n\lt\{EV^{p+1/2}_{V(n),T}-\sum^\infty_{k=1}\frac 1{2^k}\sum_{x:|x|=k} 
EV^{p+1/2}_{V(n),T_x}\rt\}<\infty. \label{e8}
\end{eqnarray}
\end{conjecture}
For an even $p\ge1$, condition \eqref{e8} becomes a hypothesis. The proof of this conjecture 
should follow the arguments of the theorem obtained in Section \ref{s1}. The 
difficulties are mainly of computational character.

The next observation is that the above-mentioned difference derivatives have depend 
more on the ``tree structure" of the operator $H_V$. Namely, the conjecture is that 
$\d V$ ($\d^kV$) have to be replaced by $\ti\d V$ ($\ti\d^kV$, respectively), where
$$
(\ti\d V)(n,j)=V(n,j)-\frac12\lt(V(n+1, 2j-1)+V(n+1,2j)\rt).
$$

The following conjecture contains the previous ones as a very particular case. This is 
a carry over of Simon's conjecture \cite[Sect.~2.8]{si2}, supplemented by 
Nazarov-Peherstorfer-Volberg-Yuditskii \cite[Lemma 6.8]{nvyu}. To formulate it, we define an 
isometry $W=\mathrm{diag}\, \{A_n\}_n: l^2\to l^2(T)$, where $A_n^t=\underbrace{[(1/2)^{n/2},\dots,(1/2)^{n/2}]}_{2^n\ \mathrm{entries}}$.

It is well-known \cite{aiz1}, that if $V$ is a radially symmetric potential, that is, 
$V(x)=const$ for a fixed $|x|$, then
$$
J_{V/\sqrt 2}=
\begin{bmatrix}
\frac{V_0}{\sqrt 2}&1&0&\ldots\\
1&\frac{V_1}{\sqrt 2}&1&\ldots \\
0&1&\frac{V_2}{\sqrt 2}&\ldots\\
\vdots&\vdots&\vdots&\ddots
\end{bmatrix}
=\frac 1{\sqrt 2} W^*H_V W,
$$
or $J_{V/\sqrt 2}$ is unitarily equivalent to a proper restriction of  $H_V$.

We rewrite the conjecture for $J_{V/\sqrt 2}$ given in \cite{nvyu} in the following 
way. Let $A(x)=\sum^{2m}_{k=0} a_k x^k$ be a polynomial of degree $2m$ and $A(x)\ge 0$ 
on $\br$. Let also
\begin{eqnarray}\label{e81}
\quad dH_V&=&\Big[V(0); \underbrace{0, 0}_{2\ \mathrm{entries}}; \underbrace{\frac12 (V(1,1), V(1,2)),\ \frac12 (V(1,1), V(1,2))}_{2\times 2=4\ \mathrm{entries}};\\
&& \underbrace{0,\dots,0}_{8\ \mathrm{entries}}\ ;
\underbrace{\frac14(V(2,1), V(2,2), V(2,3), V(2,4)),\ \dots}_{4\times 4=16\ \mathrm{entries}}\ ;\  \dots \Big]^t. \nonumber
\end{eqnarray}
Consider an operator-valued polynomial
$$
B(H_0)=\sum^{2m}_{k=0} \frac{a_k}{2^{k/2}} \underbrace{H_0W\cdot W^*H_0W\cdot \dots 
\cdot W^*H_0}_{k\mathrm{\ factors}}. 
$$
\begin{conjecture}\label{c2} Let
$$
\sum^\infty_{n=1}\frac 1{2^n}\sum_{x:|x|=n} |V(x)|^{2m+2}<\infty, \quad 
|(B(H_0)dH_V,dH_V)|<\infty.
$$
Then
\begin{eqnarray*}
&& \int^{2\sqrt 2}_{-2\sqrt 2} \log\mu'(x)\, A(x/\sqrt 2) \sqrt{8-x^2}\, dx>-\infty, 
\nonumber\\
&& \limsup_n\Big\{F^A_{V(n),T}-\sum^\infty_{k=1}\frac 1{2^k}\sum_{x:|x|=k} 
F^A_{V(n),T_x}\Big\}<\infty,
\end{eqnarray*}
where
$$
F^A_{V(n),T}=\sum_s F^A(x^+_{V(n),T; s}) + \sum_s F^A(x^-_{V(n),T; s}),
$$
and, for $\pm x^\pm>\pm 2\sqrt 2$,
\begin{equation*}
F^A(x^\pm)=\pm\int^{x^\pm}_{\pm 2\sqrt 2} A(s/\sqrt 2)\sqrt{8-s^2}\, ds.
\end{equation*}
\end{conjecture}
The notation $\{x^\pm_{V(n),T;s}\}_s$ is introduced in \eqref{e11}, \eqref{e12}.

Several remarks are in order. First, we can formulate a similar conjecture for {\it Jacobi operators} on trees.  The increments of the coefficients associated to the edges of the tree then fill in zero entries in \eqref{e81}.
Second,  consider the usual shift $Se_k=e_{k+1}$ on $l^2$, $\{e_k\}_k$ being the standard basis of the space. We now fix the standard basis $\{e_x\}_{x\in\cv(T)}$ in $l^2(T)$, see  \eqref{e01}. The binary shift is defined as
$$
S_1=
\begin{bmatrix}
0&0&0&\dots\\
B_1&0&0&\dots\\
0&B_2&0&\dots\\
\vdots&\vdots&\ddots&\ddots
\end{bmatrix},
$$
where $B_1=[1,1]^t$, and $B_{k+1}$ is the $2^{k+1}\times 2^k$ matrix $\begin{bmatrix} B_k&0\\0&B_k\end{bmatrix}$. It is plain that $H_0=S_1+{S_1}^*$ and 
\begin{eqnarray*}
S_1W&=&\sqrt2 WS,\quad W^*{S_1}^*=\sqrt2 S^*W^*,\\
W^*S_1&=&\sqrt2SW^*,\quad {S_1}^*W=\sqrt2WS^*.  
\end{eqnarray*}
For instance, we see for Theorem \ref{t1}
$$
A(x)=1,\quad B(H_0)=I,\quad (B(H_0)dH_V, dH_V)=\sum^\infty_{n=1}\frac1{2^n}\sum_{x:|x|=n} V(x)^2<\infty, 
$$
and for Theorem \ref{t2},
\begin{eqnarray*}
A(x)&=&x^2-4,\  B(H_0)=\frac 12(H_0WW^*H_0-8)=-\frac12 (S_1-{S_1}^*)^*WW^*(S_1-{S_1}^*),\\
&& (B(H_0)dH_V, dH_V)=-\sum^\infty_{n=1}\frac1{2^{n+1}}\sum_{x:|x|=n} \d V(x)^2>-\infty. 
\end{eqnarray*}
Roughly speaking, Conjecture \ref{c2} says that the role of $J_0$ is played by $H_0$ and the usual shift $S$ is replaced by $S_1$ as compared to conjectures \cite[Sect.~2.8]{si2} and \cite[Lemma 6.8]{nvyu}. The presence of the binary shift leads to the ``binary'' derivatives appearing in the formulations of the theorems.

\medskip\nt
{\it Acknowledgment.} The author would like to thank S. Denisov for helpful 
discussions.

\end{document}